\documentclass[12pt,a4paper]{article}
\usepackage[cp1251]{inputenc} 
\usepackage[russian]{babel}
\usepackage{amsfonts}

\newcommand{\di}{\displaystyle}
\newcommand{\de}{\delta}

\newcommand{\be}{\beta}
\newcommand{\ga}{\gamma}

\newcommand{\la}{\lambda}

\newcommand{\iy}{\infty}

\newcommand{\vfi}{\varphi}

\newcommand{\br}[1]{{\left\langle {#1}\right\rangle}}

\newcommand{\te}{\theta}

\begin{document}

\thispagestyle{empty}

\begin{center}
{\large\bf
SPECTRAL ANALYSIS OF THE DIRAC SYSTEM WITH A SINGULARITY IN AN INTERIOR POINT}\\[0.2cm]
{\bf O. Gorbunov, C-T. Shieh and V.Yurko} \\[0.2cm]
\end{center}

{\bf Abstract.} We study the non-selfadjoint Dirac system on the line having
an non-integrable regular singularity in an interior point with additional
matching conditions at the singular point. Special fundamental systems
of solutions are constructed with prescribed analytic and asymptotic
properties. Behavior of the corresponding Stockes multipliers is established.
These fundamental systems of solutions will be used for studying direct
and inverse problems of spectral analysis.



{\bf 1. Introduction. } Consider the Dirac system on the line with a regular singularity
at $x=0$:
$$
BY'(x)+\Big(Q_0(x)+Q(x)\Big)Y(x)=\la Y(x),\;\;-\iy<x<+\iy,               \eqno(1)
$$
where
$$
Y(x)=\left(\begin{array}{c} y_1(x) \\ y_2(x) \end{array}\right),\;
B=\left(\begin{array}{rc} 0&1\\-1&0\end{array}\right),\;
Q(x)=\left(\begin{array}{cr} q_1(x)&q_2(x)\\q_2(x)&-q_1(x)\end{array}\right),\;
Q_0(x)=\di\frac{\mu}{x}\left(\begin{array}{cc} 0&1\\1&0\end{array}\right),
$$
here $\mu$ is a complex number, $q_j(x)$ are complex-valued absolutely 
continuous functions, and $q_j'(x)\in L(-\iy,+\iy).$

In this paper special fundamental systems of solutions for system (1) are
constructed with prescribed analytic and asymptotic properties. Behavior
of the corresponding Stockes multipliers is established. These fundamental
systems of solutions will be used for studying direct and inverse problems
of spectral analysis by the contour integral method and by the method of
spectral mappings [1]-[2]. These systems can be also used for studying 
boundary value problems on a finite interval and on the half-line.

Differential equations with singularities inside the interval play an
important role in various areas of mathematics as well as in applications.
Moreover, a wide class of differential equations with turning points can be
reduced to equations with singularities. For example, such problems appear
in electronics for constructing parameters of heterogeneous electronic lines
with desirable technical characteristics [3]-[5]. Boundary value problems
with discontinuities in an interior point appear in geophysical models for
oscillations of the Earth [6]-[8]. Furthermore, direct and inverse spectral
problems for equations with singularities and turning points are used for
studying the blow-up behavior of solutions for some nonlinear integrable
evolution equations in mathematical physics (see, for example, [9]).
We also note that in different problems of natural sciences
we face different kind of matching conditions in the interior point.

The case when a singular point lies at the endpoint of the interval
was investigated fairly completely for various classes of
differential equations in [10]-[14] and other works.
The presence of singularity inside the interval produces
essential qualitative modifications in the investigation (see [15]).

A few words on the structure of the paper. In section 2 we consider
a model Dirac operator (see (2)) with the zero potential $Q(x)\equiv 0$
and without the spectral parameter. It is important that this system is
studied in the complex $x$-plane. We construct fundamental matrices for
the model system. Using analytic continuations and symmetry we calculate
directly the Stockes multipliers for the model system. Then we consider
the Dirac system on the real $x$-line with $Q(x)\equiv 0$ and with the
complex spectral parameter (see (12)), and carry over our construction
to this system. For this purpose we use a simple but important property:
if $Y(x)$ is a solution of (2), then $Y(\la x)$ is a solution of (12).
 In section 3 by perturbation theory we construct special fundamental 
matrices for system (1) with necessary analytic and asymptotic properties. 
In section 4 asymptotic properties of the Stockes multipliers
for system (1) are established.
Using these results we plan to study direct and inverse problems of
spectral analysis for system (1) in a separate paper.

\medskip
{\bf 2. Model Dirac system in the complex $x$-plane. } 
Let for definiteness, $Re\,\mu>0,$ $1/2-\mu\notin{\bf N}$
(other cases require minor modifications). Consider the model Dirac 
system without spectral parameter in the complex $x$-plane:
$$
BY'(x)+Q_0(x)Y(x)=Y(x).                                                 \eqno(2)
$$
Let $x=re^{i\vfi}, r>0, \vfi\in(-\pi,\pi],$ $x^{\xi}=
\exp(\xi(\ln r+i\vfi)),$ and $\Pi_{-}$ be the $x$-plane with the cut
$x\le 0.$ Let numbers $c_{10}, c_{20}$ be such that $c_{10}c_{20}=1.$
Then equation (2) has the matrix solution
$$
C(x)=\widehat C(x)H(x),
$$
where
$$
H(x)=\left(\begin{array}{cc}
x^{\mu_1}&0\\0&x^{\mu_2}\end{array}\right),\;
\;\widehat C(x)= \di\sum_{k=0}^{\iy}x^{2k}\left(\begin{array}{cc}
xc_{1,2k+1}&c_{2,2k}\\-c_{1,2k}&xc_{2,2k+1}\end{array}\right),
$$
$$
c_{j,2k}=(-1)^k\di\frac{c_{j0}}{2^kk!\di\prod_{s=0}^{k-1}(2\mu_j+1+2s)},\;
c_{j,2k+1}=(-1)^k\di\frac{c_{j0}}{2^kk!\di\prod_{s=0}^k(2\mu_j+1+2s)},\;
$$
$\mu_j=(-1)^j\mu,\;j=1,2.$ We agree that if a certain symbol denotes
a matrix solution of the system, then the same symbol with one index
denotes columns of the matrix, and this symbol with two indeces denotes
entries, for example, $C(x)=\Big(C_1(x),\;C_2(x)\Big)=\left(\begin{array}
{cc}C_{11}(x)&C_{12}(x)\\C_{21}(x)&C_{22}(x)\end{array}\right).$

The functions $\widehat C_k(x), k=1,2,$ are entire in $x,$ and the functions
$C_k(x), k=1,2$ are regular in $\Pi_{-}$. The functions $C_k(x), k=1,2,$
form the fundamental system of solutions for (2), and $\det C(x)\equiv 1.$

\smallskip
Denote
$$
I=\left(\begin{array}{cc}
1&0\\0&1\end{array}\right),\; J=\left(\begin{array}{cc}
0&1\\1&0\end{array}\right),\; K=\left(\begin{array}{cr}
1&0\\0&-1\end{array}\right),\;
e^0(x)=\left(\begin{array}{cr}ie^{ix}&-ie^{-ix}\\
e^{ix}&e^{-ix}\end{array}\right).
$$
Clearly, $K^2=J^2=-B^2=I$, $Q(x)=q_1(x)K+q_2(x)J,\;
Q_0(x)=\di\frac{\mu}{x}J$, $\det e^0(x)\equiv2i,$
$$
\left(\Big(I+\di\frac{a}{x}J\Big)^{-1}\right)'=
\Big(I+\di\frac{a}{x}J\Big)^{-2}\di\frac{a}{x^2}.                       \eqno(3)
$$
Note that the matrix $e^0(x)$ is a solution of the system $BY'(x)=Y(x).$

The matix Jost-type solution $e(x)$ of system (2) is constructed
from the following system of integral equations:
$$
e(x)=\Big(I-\di\frac{1}{2}Q_0(x)\Big)^{-1}e^0(x)\left(I+\di\frac{1}{2}
\int_x^\iy e^{0,-1}(x)\Big(Q'_0(t)+Q_0(t)BQ_0(t)\Big)e(t)\,dt\right),   \eqno(4)
$$
where $e^{0,-1}(t)=(e^0(t))^{-1}$. Let us show that if $e(x)$ is
a solution of equation (4), then $e(x)$ is a solution of system (2).
Denote
$D(t)=\frac{1}{2}e^{0,-1}(t)\Big(Q'_0(t)+Q_0(t)BQ_0(t)\Big)e(t).$
Then (4) takes the form $e(x)=\Big(I-\di\frac{1}{2}Q_0(x)\Big)^{-1}e^0(x)
\left(I+\int_x^\iy D(t)\,dt\right),$
and consequently,
$$
Be'(x)-e(x)=B\left(\Big(I-\di\frac{1}{2}Q_0(x)\Big)^{-1}\right)'e^0(x)
\left(I+\int_x^\iy D(t)\,dt\right)
$$
$$
+B\Big(I-\di\frac{1}{2}Q_0(x)\Big)^{-1}
\left(e^0(x)\right)'\left(I+\int_x^\iy D(t)\,dt\right)+
B\Big(I-\di\frac{1}{2}Q_0(x)\Big)^{-1}e^0(x)\Big(-D(x)\Big)-e(x).
$$
Using (3) and the relation $B\Big(I-\di\frac{1}{2}Q_0(x)\Big)^{-1}=
\Big(I+\di\frac{1}{2}Q_0(x)\Big)^{-1}B,$ we obtain
$$
Be'(x)-e(x)=\di\frac{1}{2}\Big(I+\di\frac{1}{2}Q_0(x)\Big)^{-1}BQ'_0(x)e(x)
+\Big(I+\di\frac{1}{2}Q_0(x)\Big)^{-1}\Big(I-\di\frac{1}{2}Q_0(x)\Big)e(x)
$$
$$
-\di\frac{1}{2}\Big(I+\di\frac{1}{2}Q_0(x)\Big)^{-1}B
\Big(Q'_0(x)+Q_0(x)BQ_0(x)\Big)e(x)-
\Big(I+\di\frac{1}{2}Q_0(x)\Big)^{-1}\Big(I+\di\frac{1}{2}Q_0(x)\Big)e(x)
$$
or
$$
Be'(x)-e(x)=-\Big(I+\di\frac{1}{2}Q_0(x)\Big)^{-1}
\Big(I+\di\frac{1}{2}Q_0(x)\Big)P_0(x)e(x)=-Q_0(x)e(x),
$$
i.e. $e(x)$ is a solution of system (2).

Now we go on to the solvability of equation (4). Put
$z_j(x)=e^{-R_jx}e_j(x),\;z^0_j=e^{-R_jx}e^0_j(x)=(R_j,\,1)^T,\;j=1,2,$
 here $R_1=i,$ $R_2=-i,$ and $T$ is the sign for the transposition. Since
$$
\Big(I-\di\frac{1}{2}Q_0(x)\Big)^{-1}=
\di\frac{1}{d(x)}\Big(I+\di\frac{1}{2}Q_0(x)\Big),\;
d(x):=\det\Big(I-\di\frac{1}{2}Q_0(x)\Big)=1-\di\frac{\mu^2}{4x^2},
$$
then equation (4) takes the form
$$
e(x)=\di\frac{1}{d(x)}\Big(I+\di\frac{\mu}{2x}J\Big)e^0(x)\left(
I-\frac{1}{2}\int_x^\iy e^{0,-1}(t)(J+\mu B)\frac{\mu}{t^2}e(t)\,dt\right),
$$
hence
$$
z_j(x)=\di\frac{1}{d(x)}\Big(I+\di\frac{\mu}{2x}J\Big)\left(
z^0_j-\frac{1}{2}\int_x^\iy
g^j(x,t)(J+\mu B)\frac{\mu}{t^2}z_j(t)\,dt\right),\;j=1,2,             \eqno(5)
$$
where $g^j(x,t)=e^0(x)e^{0,-1}(t)e^{R_j(t-x)},$ or
$$
g^1(x,t)=\di\frac{1}{2i}(iI-B)+\frac{1}{2i}(iI+B)e^{2i(t-x)},\;
g^2(x,t)=\di\frac{1}{2i}(iI-B)e^{-2i(t-x)}+\frac{1}{2i}(iI+B).         \eqno(6)
$$

{\bf Theorem 1. }{\it Equations (5) have analytic in $\Pi_{-}$
solutions, and }\\
{\it 1) $|z_1(x)-z^0_1|\le C/|x|$ for $|x|\ge x_0,\;
\arg x\in[-\pi+\de_0,\,\pi],$\\
2) $|z_2(x)-z^0_2|\le C/|x|$ for $|x|\ge x_0,\;
\arg x\in[-\pi,\,\pi-\de_0],$\\
where the constant $C$ depends only on $x_0,\;\de_0,\;\mu,$ and}
$x_0\sin\de_0\ge4\pi|\mu|\Big(1+|\mu|\Big).$

{\it Proof. } In view of (6), the contour in (5) for $z_1(x)$ must
be chosen such that $Im(t-x)\ge0,$ and for $z_2(x)$ such that
$Im(t-x)\le0.$ We consider two cases.

1) We choose the contour such that $\arg t=\arg x,\;|t|\ge|x|;$
then $Im(t-x)\ge0$ for $Im\,x\ge0,$ and $Im(t-x)\le0$ for $Im\,x\le0,$
i.e. $z_1(x)$ is considered for $Im\,x\ge0,$ and $z_2(x)$ --
for $Im\,x\le0.$ Denote $A(x):=(d(x))^{-1}(I+\frac{\mu}{2x}J).$
Let $x=Re^{i\te},\; t=re^{i\te},$, then (5) takes the form
$$
z_j(Re^{i\te})=A(Re^{i\te})\left(z^0_j-\frac{1}{2}
\int_R^\iy g^j(Re^{i\te},re^{i\te})(J+\mu B)
\frac{\mu e^{-i\te}}{r^2} z_j(re^{i\te})\,dr\right),\;j=1,2.            \eqno(7)
$$
We solve (7) by the method of successive approximations:
$$
\left.\begin{array}{l}
z_j(Re^{i\te})=\di\sum_{k=0}^\iy(z_j)_k(Re^{i\te}),\;\;
(z_j)_0(Re^{i\te})=A(Re^{i\te})z^0_j,\\
(z_j)_{k+1}(Re^{i\te})=-\di\frac{1}{2}A(Re^{i\te})\int_R^\iy
g^j(Re^{i\te},re^{i\te})(J+\mu B)\frac{\mu e^{-i\te}}{r^2}
(z_j)_k(re^{i\te})\,dr,\;j=1,2.
\end{array}\right\}                                                     \eqno(8)
$$
By induction we obtain $|(z_j)_k(Re^{i\te})|\le 2^{k+2}(1+|\mu|)^k/k!$
for $|x|=R\ge|\mu|.$ Therefore, the series in (8) converges uniformly
for $|x|\ge\mu$ and $Im\,x\ge0,\;Im\,x\le0$ for $z_1(x)$ and $z_2(x),$
respectively, and $z_1(x)$ is analytic for $|x|>|\mu|,\;Im\,x>0,$ and
$z_2(x)$ is analytic for $|x|>|\mu|,\;Im\,x<0;$; they are continuous
in the closure of these domains. This alows one to deform the contour
in (5) in the domain of analyticity. Moreover, one gets $|z_j(x)|\le C$
in the corresponding domain. Taking (7) into account we deduce
$$
z_j(Re^{i\te})-z^0_j=\Big(A(Re^{i\te})-I\Big)z^0_j-
\frac{1}{2}A(Re^{i\te})\int_R^\iy g^j(Re^{i\te},re^{i\te})
(J+\mu B)\frac{\mu e^{-i\te}}{r^2} z_j(re^{i\te})\,dr.
$$
Since $A(x)-I=A(x)(I-(I-\frac{\mu}{2x}J))=A(x)\frac{\mu}{2x}J,$
it follows that
$$
z_j(Re^{i\te})-z^0_j=\frac{1}{2}A(Re^{i\te})\left(
\frac{\mu}{R}e^{-i\te}Jz^0_j-\int_R^\iy g^j(Re^{i\te},re^{i\te})
(J+\mu B)\frac{\mu e^{-i\te}}{r^2} z_j(re^{i\te})\,dr\right),
$$
and consequently,
$$
|z_j(Re^{i\te})-z^0_j|\le\frac{1}{2}\cdot2\Big(\frac{|\mu|}{|R|}\cdot2+
|\mu|(1+|\mu|)\frac{C}{R}\Big)\;\;\;\mbox{or}\;\;\;
|z_j(x)-z^0_j|\le\frac{C}{|x|}.
$$

2) In (5) we take the contour $t=x+\xi,\;\xi\ge0,$ then $Im(t-x)=0,$
and (5) takes the form
$$
z_j(x)=A(x)\left(z^0_j-\frac{1}{2}\int_0^\iy
g^j(0,\xi)(J+\mu B)\frac{\mu}{(x+\xi)^2}z_j(x+\xi)\,d\xi\right).     \eqno(9)
$$
We solve (9) by the method of successive approximations:
$$
\left.\begin{array}{l}
z_j(x)=\di\sum_{k=0}^\iy(z_j)_k(x),\;\mbox{where}\;(z_j)_0(x)=A(x)z^0_j,
\\ (z_j)_{k+1}(x)=-\di\frac{1}{2}A(x)\int_0^\iy
g^j(0,\xi)(J+\mu B)\frac{\mu}{(x+\xi)^2}(z_j)_k(x+\xi)\,d\xi,\;j=1,2.
\end{array}\right\}                                                  \eqno(10)
$$
Let us prove by induction that for $(z_j)_k(x)$ from (10) for
$|x|\ge|\mu|$ one has
$$
|(z_j)_k(x)|\le4\Big(2\pi|\mu|\frac{1+|\mu|}{|x|}\Big)^k,\;Re\,x\ge0\;\;
\mbox{and}\;\;
|(z_j)_k(x)|\le4\Big(2\pi|\mu|\frac{1+|\mu|}{|Im\,x|}\Big)^k,\;Re\,x\le0.
$$
The first step is obvious. Now we assume that the estimates are valid
for $(z_j)_k(x),$ and prove them for $(z_j)_{k+1}(x).$

For $|x|\ge|\mu|,$ we have $|A(x)|\le2$ and $|g^j(0,\xi)|\le2;$ then
it follows from (10) that
$$
|(z_j)_{k+1}(x)|\le2|\mu|(1+|\mu|)\int_0^\iy\frac{1}{|x+\xi|^2}
|(z_j)_k(x+\xi)|\,d\xi,
$$
One has
$$
\int_0^\iy\frac{d\xi}{|x+\xi|^2}\le\frac{\pi}{|x|}
\;\;\mbox{for}\;\; Re\,x\ge0,\;\;\mbox{and}\;\;
\int_0^\iy\frac{d\xi}{|x+\xi|^2}\le\frac{\pi}{|Im\,x|}
\;\;\mbox{for}\;\;Re\,x\le0.                                             \eqno(11)
$$

a) For $Re\,x\ge0,$
$$
|(z_j)_{k+1}(x)|\le4\Big(2|\mu|(1+|\mu|)\Big)^{k+1}\pi^k\int_0^\iy
\frac{1}{|x+\xi|^2}\cdot\frac{1}{|x+\xi|^k}\,d\xi.
$$
Taking $\di\int_0^\iy\frac{1}{|x+\xi|^2}\cdot\frac{1}{|x+\xi|^k}\,d\xi\le
\frac{1}{|x|^k}\int_0^\iy\frac{1}{|x+\xi|^2}\,d\xi,$ and (11) into account,
we obtain our result.

\smallskip
b) For $Re\,x\le0,$
$$
|(z_j)_{k+1}(x)|\le4\Big(2|\mu|(1+|\mu|)\Big)^{k+1}\pi^k\int_0^\iy
\frac{1}{|x+\xi|^2}\cdot\frac{1}{|Im(x+\xi)|^k}\,d\xi.
$$
Since $Im(x+\xi)=Im\,x,$ one has $\di\int_0^\iy\frac{1}{|x+\xi|^2}
\cdot\frac{1}{|Im(x+\xi)|^k}\,d\xi=\frac{1}{|Im(x+\xi)|^k}\int_0^\iy
\frac{1}{|x+\xi|^2}\,d\xi.$ Using (11), we obtain our result.

Combining the results for $Re\,x\ge0$ and $Re\,x\le0,$ we can write
$$
|(z_j)_k(x)|\le4\Big(2\pi|\mu|\frac{1+|\mu|}{|x|\sin\de_0}\Big)^k,\;
|\arg x|\le\pi-\de_0.
$$
Let
$$
2\pi|\mu|\frac{1+|\mu|}{|x|\sin\de_0}\le\frac{1}{2}\;\;\mbox{or}\;\;
|x|\ge x_0=4\pi|\mu|\frac{1+|\mu|}{\sin\de_0}.
$$
Then the series (10) is majorized by the numerical convergent series.
Analogously we get $|z_j(x)-z^0_j|\le C/|x|$ for $|x|\ge x_0,\;|\arg x|
\le\pi-\de_0$. Theorem 1 is proved.

\smallskip
{\bf Corollary. }{\it $e(x)$ is a fundamental matrix, and $\det e(x)=2i.$}

\smallskip
The following lemma is important for calculating the Stockes multipliers.

\smallskip
{\bf Lemma 1. }{\it For $x\in D_+=\{z|\arg z\in(0,\pi]\}$ the following
relations hold}
$$
-Ke_2(-x)\equiv e_1(x)\;\;\;\;KC_j(-x)\equiv(-1)^je^{-i\pi\mu_j}C_j(x),\;j=1,2.
$$

{\it Proof. } Note that if $Y_0(x)$ is a solution of system (2), then
$KY_0(-x)$ is also a solution of (2). We consider integral equations
(7) for $z_1(x)$ and $z_2(-x)$ for $x\in D_+$. Let $x=Re^{i\te}.$
Then $-x=Re^{i(\te-\pi)}$):
$$
z_1(Re^{i\te})=A(Re^{i\te})\left(z^0_1-\frac{1}{2}\int_R^\iy
g^1(Re^{i\te},re^{i\te})(J+\mu B)\frac{\mu e^{-i\te}}{r^2}
z_1(re^{i\te})\,dr\right),
$$
$$
z_2(Re^{i(\te-\pi)})=A(Re^{i(\te-\pi)})\left(z^0_2-\frac{1}{2}\int_R^\iy
g^2(Re^{i(\te-\pi)},re^{i(\te-\pi)})(J+\mu B)\frac{\mu e^{-i(\te-\pi)}}{r^2}
z_2(re^{i(\te-\pi)})\,dr\right).
$$
One has
$KA(-x)=K\frac{1}{d(-x)}(I-Q_0(x)/2)$, $d(-x)=d(x)$, $KQ_0(x)=-Q_0(x)K$,
$KB=-BK$, $KA(-x)=A(x)K$, $Kg^2(-x,-t)=K(\frac{1}{2i}(iI-B)e^{2i(t-x)}+
\frac{1}{2i}(iI+B))$, then $Kg^2(-x,-t)=(\frac{1}{2i}(iI+B)e^{2i(t-x)}
+\frac{1}{2i}(iI-B))K=g^1(x,t)K$.
Multiply the second relation by $K$:
$$
Kz_2(-Re^{i\te})=A(Re^{i\te})\left(Kz^0_2-\frac{1}{2}\int_R^\iy
g^1(Re^{i\te},re^{i\te})(J+\mu B)\frac{\mu e^{-i\te}}{r^2}
Kz_2(-re^{i\te})\,dr\right).
$$
Since $Kz^0_2=-z^0_1$, then for the function $\widetilde z_2(x)=-Kz_2(-x),$
we have the relation
$$
\widetilde z_2(Re^{i\te})=A(Re^{i\te})\left(z^0_1-\frac{1}{2}\int_R^\iy
g^1(Re^{i\te},re^{i\te})(J+\mu B)\frac{\mu e^{-i\te}}{r^2}
\widetilde z_2(re^{i\te})\,dr\right).
$$
The functions $\widetilde z_2(x)$ and $z_1(x)$ satisfy the same equation;
This yeilds $\widetilde z_2(x)\equiv z_1(x).$ Taking the relation
$e_j(x)=e^{R_jx}z_j(x)$ into account, we obtain the first assertion of the lemma.

Furthermore, since $C_j(x)=x^{\mu_j}\widehat C_j(x),$ it follows that
$C_j(-x)=(-x)^{\mu_j}\widehat C_j(-x).$ Moreover,\\
$-x=\left\{\begin{array}{ll} xe^{i\pi}&\mbox{ for }\arg x\in(-\pi,0],\\
xe^{-i\pi}&\mbox{ for }\arg x\in(0,\pi].\end{array}\right.$
This yields $(-x)^{\mu_j}=\left\{\begin{array}{ll}
x^{\mu_j}e^{i\pi\mu_j}&\mbox{ for }\arg x\in(-\pi,0],\\
x^{\mu_j}e^{-i\pi\mu_j}&\mbox{ for }\arg x\in(0,\pi].\end{array}\right.$

Thus, for $x\in D_+$ one has $C_j(-x)=e^{-i\pi\mu_j}x^{\mu_j}
\widehat C_j(-x).$ Then $K\widehat C_j(-x)=(-1)^j\widehat C_j(x),$
and the lemma is proved.

\smallskip
In the domain $|\arg x|\le\pi-\de_0$ we have two fundamental matrices; then
$e(x)=C(x)\ga^0$ and $C(x)=e(x)\be^0$; the matrices $\ga^0,\;\be^0$ are called
the Stockes multipliers.

\smallskip
{\bf Theorem 2. }{\it For the Stockes multipliers of system (2) the following
relations hold} $\det\ga^0=2i$, $\ga^0_{11}=e^{-i\pi\mu_1}\ga^0_{12},\;
\ga^0_{21}=-e^{-i\pi\mu_2}\ga^0_{22},\;\ga^0_{11}\ga^0_{21}=(i\cos\pi\mu)^{-1}.$

\smallskip
{\it Proof. }
The first assertion follows from the relations $\det e(x)=\det C(x)\det\ga^0,$
$\det e(x)\equiv2i,$ $\det C(x)\equiv1.$ In order to prove the second assertion
we rewrite $e(x)=C(x)\ga^0$ in the vector form:
$$
e_1(x)=\ga^0_{11}C_1(x)+\ga^0_{21}C_2(x),\;e_2(x)=\ga^0_{12}C_1(x)+\ga^0_{22}C_2(x).
$$
Let $x\in D_+$. Substututing $-x$ to the second relation and multiplying on
$(-K),$ we get
$e_1(x)=\ga^0_{12}e^{-i\pi\mu_1}C_1(x)+\ga^0_{22}(-e^{-i\pi\mu_2})C_2(x).$
Therefore $\ga^0_{11}=e^{-i\pi\mu_1}\ga^0_{12},\;
\ga^0_{21}=-e^{-i\pi\mu_2}\ga^0_{22}$. Since $\det\ga^0=\ga^0_{11}
\cdot (-e^{i\pi\mu_2})\ga^0_{21}-e^{i\pi\mu_1}\ga^0_{11}\ga^0_{21},$
it follows thst $\ga^0_{11}\ga^0_{21}=(i\cos\pi\mu)^{-1}$. Theorem 2 is proved.

\smallskip
{\bf Corollary. }{\it The following properties of the Stockes multipliers
$\be^0$ hold:\\
$\det\be^0=(2i)^{-1}$, $\be^0_{11}=e^{-i\pi\mu_1}\be^0_{21},\;\be^0_{12}=-e^{-i\pi\mu_2}\be^0_{22},
\;\be^0_{21}\be^0_{22}=(4i\cos\pi\mu)^{-1}.$}

\smallskip
Now we consider the sytem
$$
BY'+Q_0(x)Y=\la Y.                                                          \eqno(12)
$$
for real $x\not=0$ and complex $\la.$ We will use a simple but important
property: if $Y(x)$ is a solution of (2), then $Y(\la x)$ is a solution of (12).

Denote $C(x,\la)=C(x\la)H(\la^{-1})$, $e(x,\la)=e(x\la).$ Clearly,
$C_j(x,\la)=x^{\mu_j}\widehat C_j(x,\la),$ where $\widehat C_j(x,\la)=
\widehat C_j(x\la)$, $e_j(x,\la)=e^{R_j\la x}z_j(x\la),\; j=1,2.$
The following theorem is obvious.

\smallskip
{\bf Theorem 3. }{\it
1) $C(x,\la)$ is a fundamental matrix for system (12), $\det C(x,\la)\equiv1,$
$C(x,\la)$ is entire in $\la,$ and $|\widehat C(x\la)|\le C$ for each $x\la$
from a compact.

2) $e(x,\la)$ is a fundamental matrix for system (12), $\det e(x,\la)\equiv2i,$
and \\ $|z_j(x\la)-z^0_j|\le C_0|x\la|^{-1}$ for $|x\la|\ge
x_0,\;\arg(x\la)\in [-\pi+\de_0,\pi]$ for $j=1,$
$\arg(x\la)\in[-\pi,\pi-\de_0]$ for $j=2,$ where $C_0$ depends only on
$x_0,\;\mu,\;\de_0$, and $x_0\sin\de_0\ge4\pi|\mu| (1+|\mu|).$

3) Let $e(x,\la)=C(x,\la)\ga^0(\la)$ and $C(x,\la)=e(x,\la)\be^0(\la).$ Then\\
$\ga^0_{jk}(\la)=\la^{\mu_j}\ga^0_{jk},\;\be^0_{kj}(\la)=\la^{-\mu_j}\be^0_{kj},
\;k,j=1,2.$}

\bigskip
{\bf 3. Fundamental systems of solutions. } Now we consider system (1) and assume
that  $\di\int_{|x|\le1}|x|^{-2Re\mu}|Q(x)|\,dx+\int_{|x|\ge1}|Q(x)|\,dx<\iy$.
 In this section we construct fundamental matrices for system (1) and 
establish properties of their Stockes multipliers. The following assertion 
is proved by the well-known method (see, for example, [1]-[2]).

\smallskip
{\bf Theorem 4. }{\it System (1) has a fundamental system of solutions
$S_j(x,\la)=x^{\mu_j}\widehat S_j(x,\la),\;j=1,2,$ where the functions
$\widehat S_j(x,\la)$ are solutions of the integral Volterra equations (13):
$$
\widehat S_j(x,\la)=\widehat C_j(x,\la)+\int_0^x C(x,\la)C^{-1}(t,\la)
\Big(\di\frac{t}{x}\Big)^{\mu_j}BQ(t)\widehat S_j(t,\la)\,dt,\;j=1,2.     \eqno(13)
$$
The functions $S_j(x,\la)$ are entire in $\la,$ and
$|\widehat S_j(x,\la)|\le C$ on compacts.}

\smallskip
Let us now construct the Birkhoff-type fundamental system of solutions
for system (1). For definiteness, we confine ourselves to the case $x>0.$
In section 2 we constructed the solution $e(x,\la)$ of equation (12) for
$|x\la|\ge x_0,\;|\arg\la|\le\pi-\de_0$, where $x_0>0,\;\de_0>0$ are such
that $x_0\sin\de_0\ge4\pi |\mu|(1+|\mu|).$ The Stockes multipliers allow
one to extend this solution by $e(x,\la)=C(x,\la)\ga^0(\la)$ on $\Pi_{-}$
and $x\not=0.$
Denote \\
$F(x\la)=\left(\begin{array}{cc}F_1(x\la)&0\\0&F_2(x\la)
\end{array}\right)$, $F_j(x\la)=\left\{\begin{array}{cl}(x\la)^{-\mu}&
\mbox{for}\;|x\la|<2|\mu|,\\e^{R_j\la x},&\mbox{for}\;|x\la|\ge2|\mu|,
\end{array}\right.,\;R_1=i,\;R_2=-i.$
Let $U^0(x,\la)=(U^0_1(x,\la),\;U^0_2(x,\la)):=e(x,\la)F^{-1}(x\la).$
It is easy to check that $|U^0(x,\la)|\le C$ for $x>0,\;|\arg\la|\le\pi/2.$
The Birkhoff-type solutions $E_j(x,\la),\;j=1,2,$ of system (1) is constructed
from the following systems of integral equations: \\
1) for $x\le a_\la:=2|\mu|/|\la|$
$$
E_1(x,\la)=e_1(x,\la)+e(x,\la)\Big(I_1\int_0^xe^{-1}(t,\la)BQ(t)E_1(t,\la)\,dt
-I_2\int_x^{a_\la}e^{-1}(t,\la)BQ(t)E_1(t,\la)\,dt
$$
$$
-\frac{1}{2}I_2e^{-1}(a_\la,\la)Q^{-1}(a_\la,\la)Q(a_\la)E_1(a_\la,\la)\Big),  \eqno(14)
$$
$$
E_2(x,\la)=e_2(x,\la)+e(x,\la)\di\int_0^xe^{-1}(t,\la)BQ(t)E_2(t,\la)\,dt;     \eqno(15)
$$
2) for $x\ge a_\la$
$$
E_1(x,\la)=e_1(x,\la)-\frac{1}{2}Q^{-1}(x,\la)Q(x)E_1(x,\la)
$$
$$
+e(x,\la)\Big(I_1\int_0^{a_\la}e^{-1}(t,\la)BQ(t)E_1(t,\la)\,dt
+\frac{1}{2}I_1\int_{a_\la}^xe^{-1}(t,\la)L(t,\la)E_1(t,\la)\,dt
$$
$$
-\frac{1}{2}I_2\int_x^\iy e^{-1}(t,\la)L(t,\la)E_1(t,\la)\,dt
+\frac{1}{2}I_1e^{-1}(a_\la,\la)Q^{-1}(a_\la,\la)Q(a_\la)E_1(a_\la,\la)\Big), \eqno(16)
$$
$$
E_2(x,\la)=e_2(x,\la)-\frac{1}{2}Q^{-1}(x,\la)Q(x)E_2(x,\la)
+e(x,\la)\Big(\di\int_0^{a_\la}e^{-1}(t,\la)BQ(t)E_2(t,\la)\,dt
$$
$$
+\frac{1}{2}\int_{a_\la}^xe^{-1}(t,\la)L(t,\la)E_2(t,\la)\,dt
+\frac{1}{2}e^{-1}(a_\la,\la)Q^{-1}(a_\la,\la)Q(a_\la)E_2(a_\la,\la)\Big),   \eqno(17)
$$
where $I_1=\left(\begin{array}{rc} 1&0\\0&0\end{array}\right),\;
I_2=\left(\begin{array}{rc} 0&0\\0&1\end{array}\right),\;Q(x,\la)=Q_0(x)-\la I,$
$$
L(t,\la)=\Big(Q^{-1}(t,\la)Q(t)\Big)'+Q^{-1}(t,\la)
\Big(Q(t)BQ(t)+Q(t)BQ(t,\la)+Q(t,\la)BQ(t)\Big).                            \eqno(18)
$$

Let us show that if $E_j(x,\la),\,j=12$ are solutions of these systems,
then they are solutions of (1). Since $Be'(x,\la)+Q(x,\la)e(x,\la)=0,$
it follows from (14)-(15) that for $x\le a_\la,$
$$
BE_j'(x,\la)+Q(x,\la)E_j(x,\la)=B(BP(x)E_j(x,\la)).
$$
Together with $B^2=-I$ this yields that for $x\le a_\la$ the functions
$E_j(x,\la)$ are solutions of system (1).

For $x\ge a_\la,$ it follows from (16)-(17) that
$$
BE_j'(x,\la)+Q(x,\la)E_j(x,\la)=-\frac{1}{2}B\Big(Q^{-1}(x,\la)Q(x)E_j(x,\la)
\Big)'
$$
$$
+\frac{1}{2}BL(x,\la)E_j(x,\la)-\frac{1}{2}BQ(x)E_j(x,\la).
$$
In view of (18) this yields
$$
BE_j'(x,\la)+Q(x,\la)E_j(x,\la)=-\frac{1}{2}BQ^{-1}(x,\la)Q(x)E'_j(x,\la)
+\Big(-\frac{1}{2}B\Big(Q^{-1}(x,\la)Q(x)\Big)'
$$
$$
+\frac{1}{2}B\Big(Q^{-1}(x,\la)Q(x)\Big)'+\frac{1}{2}B^2Q(x)
+\frac{1}{2}BQ^{-1}(x,\la)Q(x)B\Big(Q(x)+Q(x,\la)\Big)-\frac{1}{2}Q(x)
\Big)E_j(x,\la),
$$
or
$$
\Big(I-\frac{1}{2}BQ^{-1}(x,\la)Q(x)B\Big)\Big(BE'_j(x,\la)+
\Big(Q(x)+Q(x,\la)\Big)E_j(x,\la)\Big)=0.
$$
Thus, the functions $E_j(x,\la)$ satisfy (1) in the points $(x,\la)$
where $\det(I-\frac{1}{2}BQ^{-1}(x,\la)Q(x)B)\not=0.$ Let us show
that for $\la$ sufficiently large, this determinant differs from zero for
each $x\ge a_\la.$  Denote $d(x,\la)=\mu^2/x^2-\la^2,$ then $Q^{-1}(x,\la)
=(d(x,\la))^{-1}(\frac{\mu}{x}J+\la I).$ Using anticommutativity of the
matrices $J,\;K,\;B,$ we obtain
$$
\di\frac{1}{2}BQ^{-1}(x,\la)Q(x)B=-\di\frac{1}{2d(x,\la)}
\Big(\frac{\mu}{x}J-\la I\Big)\Big(q_1(x)K+q_2(x)J\Big)
$$
Since $J^2=I$ and $JK=-B,$ it follows that
$$
\det\Big(I-\frac{1}{2}BQ^{-1}(x,\la)Q(x)B\Big)
$$
$$
=\det\Big(I+\di\frac{1}{2d(x,\la)}\Big(-q_1(x)\frac{\mu}{x}B+
q_2(x)\frac{\mu}{x}I-q_1(x)\la K-q_2(x)\la J\Big)\Big)
$$
$$
=\frac{1}{4d^2(x,\la)}
\left|\begin{array}{cc}
\di2d(x,\la)+q_2(x)\frac{\mu}{x}-q_1(x)\la&\di -q_1(x)\frac{\mu}{x}-q_2(x)\la\\
\di q_1(x)\frac{\mu}{x}-q_2(x)\la&\di 2d(x,\la)+q_2(x)\frac{\mu}{x}+q_1(x)\la
\end{array}\right|
$$
$$
=\frac{1}{4d^2(x,\la)}\Big(\Big(2d(x,\la)+q_2(x)\frac{\mu}{x}\Big)^2-
q_1^2(x)\la^2+q_1^2(x)\frac{\mu^2}{x^2}-q_2^2(x)\la^2\Big)
$$
$$
=\frac{1}{4d^2(x,\la)} \Big(4d^2(x,\la)+4d(x,\la)q_2(x)\frac{\mu}{x}+
\Big(q_1^2(x)+q_2^2(x)\Big)\Big(\frac{\mu^2}{x^2}-\la^2\Big)\Big),
$$
i.e.
$$
\det\Big(I-\frac{1}{2}BQ^{-1}(x,\la)Q(x)B\Big)=1+\frac{1}{4d(x,\la)}
\Big(4q(x)\frac{\mu}{x}+q_1^2(x)+q_2^2(x)\Big).
$$
We estimate the second term. For $x\ge a_\la$ we have $|d(x,\la)|\ge|\la|^2
-|\mu/x|^2\ge|\la|^2/2.$ Since $q_1(x)$ and $q_2(x)$ are bounded it follows
that
$$
\left|\det\Big(I-\frac{1}{2}BQ^{-1}(x,\la)Q(x)B\Big)-1\right|\le
\frac{1}{2|\la|^2}\Big(4C\frac{|\la|}{2}+2C^2\Big)\le\frac{C_0}{|\la|}.
$$
For $|\la|\ge2C_0$ we get $\det(I-\frac{1}{2}BQ^{-1}(x,\la)Q(x)B)\ge 1/2.$
Therefore, for $x\ge a_\la$ and sufficiently large $|\la|,$ the function
$E_j(x,\la)$ is a solution of system (1).

Let us go on to the solvability of systems (14)-(17). Denote
\\ $U(x,\la)=(U_1(x,\la),\;U_2(x,\la)):=E(x,\la)F^{-1}(x\la),$
where $E(x,\la)=(E_1(x,\la),\;E_2(x,\la)).$ Then for $U_j(x,\la),\;j=1,2,$
the following relations hold: 1) for $x\le a_\la,$
$$
U_1(x,\la)=U^0_1(x,\la)+e(x,\la)\Big(I_1\int_0^xe^{-1}(t,\la)BQ(t)
\frac{F_1(t\la)}{F_1(x\la)}U_1(t,\la)\,dt
$$
$$
-I_2\int_x^{a_\la}e^{-1}(t,\la)BQ(t)\frac{F_1(t\la)}{F_1(x\la)}U_1(t,\la)\,dt
-\frac{1}{2}I_2\int_{a_\la}^\iy e^{-1}(t,\la)L(t,\la)
\frac{F_1(t\la)}{F_1(x\la)}U_1(t,\la)\,dt
$$
$$
-\frac{1}{2}I_2e^{-1}(a_\la,\la)Q^{-1}(a_\la,\la)Q(a_\la)
\frac{F_1(a_\la\la)}{F_1(x\la)}U_1(a_\la,\la)\Big),                        \eqno(19)
$$
$$
U_2(x,\la)=U^0_2(x,\la)+e(x,\la)\di\int_0^xe^{-1}(t,\la)BQ(t)
\frac{F_2(t\la)}{F_2(x\la)}U_2(t,\la)\,dt;                                 \eqno(20)
$$
2) for $x\ge a_\la,$
$$
U_1(x,\la)=U^0_1(x,\la)-\frac{1}{2}Q^{-1}(x,\la)Q(x)U_1(x,\la)+
e(x,\la)\Big(I_1\int_0^{a_\la}e^{-1}(t,\la)BQ(t)\frac{F_1(t\la)}
{F_1(x\la)}U_1(t,\la)\,dt
$$
$$
+\frac{1}{2}I_1\int_{a_\la}^xe^{-1}(t,\la)L(t,\la)
\frac{F_1(t\la)}{F_1(x\la)}U_1(t,\la)\,dt
-\frac{1}{2}I_2\int_x^\iy e^{-1}(t,\la)L(t,\la)
\frac{F_1(t\la)}{F_1(x\la)}U_1(t,\la)\,dt
$$
$$
+\frac{1}{2}I_1e^{-1}(a_\la,\la)Q^{-1}(a_\la,\la)Q(a_\la)
\frac{F_1(a_\la\la)}{F_1(x\la)}U_1(a_\la,\la)\Big),                       \eqno(21)
$$
$$
U_2(x,\la)=U^0_2(x,\la)-\frac{1}{2}Q^{-1}(x,\la)Q(x)U_2(x,\la)
+e(x,\la)\Big(\di\int_0^{a_\la}e^{-1}(t,\la)Q(t)
\frac{F_2(t\la)}{F_2(x\la)}U_2(t,\la)\,dt
$$
$$
+\frac{1}{2}\int\limits_{a_\la}^xe^{-1}(t,\la)L(t,\la)
\frac{F_2(t\la)}{F_2(x\la)}U_2(t,\la)\,dt
+\frac{1}{2}e^{-1}(a_\la,\la)Q^{-1}(a_\la,\la)Q(a_\la)
\frac{F_2(a_\la\la)}{F_2(x\la)}U_2(a_\la,\la)\Big).                     \eqno(22)
$$
Since $e^{-1}(x,\la)=-\di\frac{1}{2i}Be^T(x,\la)B,$ it follows that
$$
e(x,\la)I_je^{-1}(t,\la)=-\di\frac{1}{2i}U^0(x,\la)
F(x\la)I_jBF^T(t\la)U^{0,T}(t,\la)B,\;j=1,2,
$$
where $U^{0,T}(t,\la)=(U^0(t,\la))^T.$ Denote $B_1=I_1B.$ Then
$F(x\la)B_1F(t\la)=F_1(x\la)F_2(t\la)B_1$. Analogously, one gets
$F(x\la)I_2BF(t\la)=F_1(t\la)F_2(x\la)B_2$, where $B_2=I_2B.$

Denote $N(x,t,\la)=F(x\la)BF(t\la)\di\frac{F_2(t\la)}{F_2(x\la)}.$ Then
$$
N(x,t,\la)=\Big(F_2(t\la)\Big)^2
\frac{F_1(x\la)}{F_2(x\la)}B_1+F_1(t\la)F_2(t\la)B_2.                  \eqno(23)
$$
We note that for $x<a_\la$ one has $F_1(x\la)=F_2(x\la).$
We rewrite (19)-(22) in the form: 1) for $x\ge a_\la,$
$$
U_1(x,\la)=U^0_1(x,\la)-\frac{1}{2}Q^{-1}(x,\la)Q(x)U_1(x,\la)
$$
$$
+\frac{1}{2i}U^0(x,\la)\Big(B_1\int_0^{a_\la}
F_1(t\la)F_2(t\la)U^{0,T}(t,\la)Q(t)U_1(t,\la)\,dt
$$
$$
-\frac{1}{2}B_1\int_{a_\la}^x
F_1(t\la)F_2(t\la)U^{0,T}(t,\la)BL(t,\la)U_1(t,\la)\,dt
$$
$$
+\frac{1}{2}B_2\int_x^\iy F^2_1(t\la)\frac{F_2(x\la)}
{F_1(x\la)}U^{0,T}(t,\la)BL(t,\la)U_1(t,\la)\,dt
$$
$$
-\frac{1}{2}B_1F_1(a_\la\la)F_2(a_\la\la)
U^{0,T}(a_\la,\la)BQ^{-1}(a_\la,\la)Q(a_\la)U_1(a_\la,\la)\Big),      \eqno(24)
$$
$$
U_2(x,\la)=U^0_2(x,\la)-\frac{1}{2}Q^{-1}(x,\la)Q(x)U_2(x,\la)
$$
$$
+\frac{1}{2i}U^0(x,\la)\Big(\di\int_0^{a_\la}
N(x,t,\la)U^{0,T}(t,\la)Q(t)U_2(t,\la)\,dt
$$
$$
-\frac{1}{2}\int_{a_\la}^xN(x,t,\la)U^{0,T}(t,\la)BL(t,\la)U_2(t,\la)\,dt
$$
$$
-\frac{1}{2}N(x,a_\la,\la)U^{0,T}(a_\la,\la)
BQ^{-1}(a_\la,\la)Q(a_\la)U_2(a_\la,\la)\Big);                      \eqno(25)
$$
2) for $x<a_\la,$
$$
U_1(x,\la)=U^0_1(x,\la)+\frac{1}{2i}U^0(x,\la)\Big(B_1\int_0^x
F_1(t\la)F_2(t\la)U^{0,T}(t,\la)Q(t)U_1(t,\la)\,dt
$$
$$
-B_2\int_x^{a_\la}F^2_2(t\la)U^{0,T}(t,\la)Q(t)U_1(t,\la)\,dt
+\frac{1}{2}B_2\int_{a_\la}^\iy F^2_2(t\la)U^{0,T}(t,\la)BL(t,\la)
U_1(t,\la)\,dt
$$
$$
+\frac{1}{2}B_2F^2_2(a_\la\la)
U^{0,T}(a_\la,\la)BQ^{-1}(a_\la,\la)Q(a_\la)U_1(a_\la,\la)\Big),     \eqno(26)
$$
$$
U_2(x,\la)=U^0_2(x,\la)+\di\frac{1}{2i}U^0(x,\la)\int_0^x
N(x,t,\la)U^{0,T}(t,\la)Q(t)U_2(t,\la)\,dt.                          \eqno(27)
$$

{\bf Lemma 3. }{\it The following estimates hold:\\
1) for $t\ge2a_\la:$\\
$|L(t,\la)|\le\di\frac{2}{|\la|}|P'(x)|+
C\Big(\frac{1}{|\la|}+\frac{t^{-\nu}}{|\la|^\nu}\Big)|P(t)|$, where
$\nu=\min\{1,\,2Re\mu\};$\\
2) for $t\le x<a_\la:\;\;\;N(x,t,\la)=(t\la)^{-2\mu}B,$

for $t<a_\la\le x:\;\;\;|N(x,t,\la)|\le|\la|^{-2Re\mu}t^{-2Re\mu},$

for $a_\la\le t\le x:\;\;\;|N(x,t,\la)|\le1.$

Proof. } Since $(Q^{-1}(x,\la))'=Q^{-2}(x,\la)\di\frac{\mu}{x^2}J,$
it follows that
$$
L(t,\la)=Q^{-1}(t,\la)\Big(Q^{-1}(t,\la)\di\frac{\mu}{t^2}JQ(t)+Q'(t)
+Q(t)BQ(t)+Q(t,\la)BQ(t)+Q(t)BQ(t,\la)\Big).
$$
It is easy to check that $KBJ=-JBK,$ and consequently,
$$
Q(t,\la)BQ(t)+Q(t)BQ(t,\la)
$$
$$
=\Big(\frac{\mu}{t}J-\la I\Big)B\Big(q_1(t)K+q_2(t)J\Big) +\Big(q_1(t)K+
q_2(t)J\Big)B\Big(\frac{\mu}{t}J-\la I\Big)=-2q_2(t)\frac{\mu}{t}B.
$$
Similarly, one gets
$$
Q(t)BQ(t)=\Big(q_1(t)K+q_2(t)J\Big)B\Big(q_1(t)K+q_2(t)J\Big)=
-\Big(q_1(t)^2+q_2(t)^2\Big)B.
$$
Substituting these relations into $L(t,\la),$ we calculate
$$
L(t,\la)=Q^{-1}(t,\la)\Big(Q'(t)-\Big(q_1(t)^2+q_2(t)^2\Big)B\Big)+
Q^{-1}(t,\la)\frac{\mu}{t}\Big(Q^{-1}(t,\la)\frac{1}{t}JQ(t)-2q_2(t)B\Big).
$$
For $t\ge a_\la$ we have
$$
|Q^{-1}(t,\la)|\le
\di\frac{|\la|+|\mu t^{-1}|}{|\la|^2-|\mu t^{-1}|^2}\le\frac{2}{|\la|}.  \eqno(28)
$$
Since $q_1(x)$ and $q_2(x)$ are bounded, it follows that
$$
|L(t,\la)|\le\frac{2}{|\la|}|Q'(t)|+
C\Big(\frac{1}{|\la|}+\frac{1}{|\la|t}\Big)|Q(t)|.
$$
If $Re\mu\ge 1/2,$ then $\nu=1,$ and our estimate is obtained;
if $0<Re\mu<1/2,$ then $\nu=2Re\mu.$ Since $1/t=t^{-\nu}t^{\nu-1}$ and
$\nu-1<0,$ it follows that $1/t\le t^{-\nu}|2\mu/\la|^{\nu-1},$ and
our estimwte is obtained too.

In order to prove the second assertion, we use (23).

a) Let $t\le x<a_\la.$ Then $F_j(t\la)=(t\la)^{-\mu},\;
F_j(x\la)=(x\la)^{-\mu},$ hence
$$
N(x,t,\la)=(t\la)^{-2\mu}B_1+(t\la)^{-2\mu}B_2.
$$

b) Let $t<a_\la\le x.$ Then $F_j(t\la)=(t\la)^{-\mu},\;
F_j(x\la)=e^{R_j\la x},$ hence
$$
N(x,t,\la)=(t\la)^{-2\mu}e^{2i\la x}B_1+(t\la)^{-2\mu}B_2.
$$
Since $x>0$ and $Im\la\ge0,$ then $|e^{2i\la x}|\le1,$ and
$|N(x,t,\la)|\le|\la t|^{-2Re\mu}.$

c) Let $a_\la\le t\le x.$ Then $F_j(t\la)=e^{R_j\la t},\;
F_j(x\la)=e^{R_j\la x},$ hence
$$
N(x,t,\la)=e^{2i\la(x-t)}B_1+B_2.
$$
Since $x-t\ge0$ and $Im\la\ge0,$ it follows that $|N(x,t,\la)|\le1.$
The lemma is proved.

\smallskip
Now we formulate and prove the main result of this section.

\smallskip
{\bf Theorem 5. }{\it Systems (24)-(25) and (26)-(27) have solutions
$U_j(x,\la),\;j=1,2$ for $x>0$ and
$\la\in\{\la:\,|\la|\ge\la_0,\;\arg\la\in(0,\pi/2]\},$ and
$|U_j(x,\la)-U^0_j(x,\la)|\le M/|\la|^\nu,$ where the constant $M$
depends on $\mu,\;Q(x),\;Q'(x).$}

\smallskip
{\bf I.} We begin with (25), (27) for $U_2(x,\la).$

a) Let $x\le a_\la.$ We construct the solution $U_2(x,\la)$ by the
method of successive approximations:
$$
U_2(x,\la)=\sum_{k=0}^\iy(U_2)_k(x,\la),\;\;\mbox{where }\;
(U_2)_0(x,\la)=U^0_2(x,\la),
$$
$$
(U_2)_{k+1}(x,\la)=\di\frac{1}{2i}U^0(x,\la)\int_0^x
N(x,t,\la)U^{0,T}(t,\la)Q(t)(U_2)_k(t,\la)\,dt.
$$
Using Lemma 3, by induction we get
$$
|(U_2)_k(x,\la)|\le\frac{C}{k!}\Big(\frac{C^2}{2|\la|^{2Re\mu}}
\int_0^{a_\la}t^{-2Re\mu}|Q(t)|\,dt\Big)^k.
$$
This means that the series converges uniformly, and consequently,
the function $U_2(x,\la)$ is continuous with respect to $x$ and
analytic with respect to $\la,$ and $|U_2(x,\la)|<C.$
Furthermore,
$$
U_2(x,\la)-U^0_2(x,\la)=\di\frac{1}{2i}U^0(x,\la)\int_0^x
N(x,t,\la)U^{0,T}(t,\la)P(t)U_2(t,\la)\,dt.
$$
Using Lemma 3, we obtain for $x\le a_\la$:
$|U_2(x,\la)-U^0_2(x,\la)|\le C/|\la|^{2Re\mu}.$

b) Let $x>a_\la.$ The solution is also found
by the method of successive approximations:
$$
U_2(x,\la)=\sum_{k=0}^\iy(U_2)_k(x,\la),\;\;\mbox{where }
$$
$$
(U_2)_0(x,\la)=U^0_2(x,\la)-\di\frac{1}{2i}U^0(x,\la)N(x,a_\la,\la)
U^{0,T}(a_\la,\la)BQ^{-1}(a_\la,\la)Q(a_\la)U_2(a_\la,\la)
$$
$$
+\di\frac{1}{2i}U^0(x,\la)\int_0^{a_\la}
N(x,t,\la)U^{0,T}(t,\la)Q(t)U_2(t,\la)\,dt,
$$
$$
(U_2)_{k+1}(x,\la)=-\frac{1}{2}Q^{-1}(x,\la)Q(x)(U_2)_k(x,\la)
$$
$$
-\frac{1}{4i}U^0(x,\la)\int_{a_\la}^x
N(x,t,\la)U^{0,T}(t,\la)BL(t,\la)(U_2)_k(t,\la)\,dt.
$$
Using results from the case a), Lemma 3 and (28), we obtain the estimates
$$
|(U_2)_0(x,\la)|\le C\Big(1+\di\frac{1}{|\la|^\nu}\Big),
$$
$$
|(U_2)_k(x,\la)|\le C\Big(1+\frac{1}{|\la|^\nu}\Big)C^k\Big(\frac{1}{|\la|}+
\frac{1}{|\la|}\int_0^\iy\Big(|Q'(t)|+|Q(t)|\Big)\,dt+
\frac{1}{|\la|^\nu}\int_0^\iy t^{-\nu}|Q(t)|)\,dt\Big)^k.
$$
For sufficiently large $|\la|\ge\la_0$, the series for $U_2(x,\la)$ converges
uniformly, hence $U_2(x,\la)$ is continuous with respect to $x$ and analytic
with respect to $\la,$ and $|U_2(x,\la)|\le C.$ Together with Lemma 3 and (28)
this yields
$$
|U_2(x,\la)-U^0_2(x,\la)|\le C\Big(\frac{1}{|\la|^{2Re\mu}}+\frac{1}{|\la|}\Big),
$$
and we arrive at the required estimate.

\smallskip
{\bf II.} Now we consider the existence of the solution $U_1(x,\la)$
of system (24), (26). The system for $U_1(x,\la)$ has the form
$$
U_1(x,\la)=U^0_1(x,\la)+D_1(x,\la)U_1(x,\la)+D_2(x,\la)U_1(a_\la,\la)+
\int_0^{+\iy}D_3(x,t,\la)U_1(t,\la)\,dt.
$$
We solve this system by the method of successive approximations:
$$
U_1(x,\la)=\sum_{k=0}^\iy(U_1)_k(x,\la),\;\;\mbox{��� }
(U_1)_0(x,\la)=U^0_1(x,\la)
$$
$$
(U_1)_{k+1}(x,\la)=D_1(x,\la)(U_1)_k(x,\la)+D_2(x,\la)(U_1)_k(a_\la,\la)+
\int_0^{+\iy}D_3(x,t,\la)(U_1)_k(t,\la)\,dt.
$$
It is easy to check that if for all $x$ the following estimates
$$
|U^0_1(x,\la)|\le D_0,\;|D_1(x,\la)|\le D_1(\la),\;|D_2(x,\la)|\le
D_2(\la),\; |D_3(x,t,\la)|\le D_3(t,\la),
$$
are vald, the
$$
|(U_1)_k(x,\la)|\le
D_0\Big(D_1(\la)+D_2(\la)+\int_0^{+\iy}D_3(t,\la)\,dt\Big)^k.           \eqno(29)
$$
Let us obtain the required estimates for the system for $U_1(x,\la).$

1) Since $|U^0_1(x,\la)|\le C,$ it follows that $D_0=C.$

2) According to the integral equation $D_1(x,\la)=0$ for $x\le a_\la,$ and\\
$D_1(x,\la)=-\frac{1}{2}Q^{-1}(x,\la)Q(x)$ for $x>a_\la.$ By virtue of (28),
$|D_1(x,\la)|\le |Q(x)|/|\la|,$ i.e. $D_1(\la)=C/|\la|.$

3) Since
$$
D_2(x,\la)=\left\{\begin{array}{ll}\di\frac{1}{4i}U^0(x,\la)B_2F^2_1(a_\la\la)
U^{0,T}(a_\la,\la)BQ^{-1}(a_\la,\la)Q(a_\la)&
\mbox{for }\; x<a_\la,\\[3mm]
\di\frac{1}{4i}U^0(x,\la)B_1F_1(a_\la\la)F_2(a_\la\la)U^{0,T}(a_\la,\la)B
Q^{-1}(a_\la,\la)Q(a_\la) &\mbox{for }\; x\ge a_\la,
\end{array}\right.
$$
it follows that $|D_2(x,\la)|\le C|Q(a_\la)|/|\la|,$ i.e. $D_2(\la)=C/|\la|.$

4) The function $D_3(x,t,\la)$ has a more complicated structure; it is
convenient to consider two cases.

a) Let $x<a_\la.$ Then
$$
D_3(x,t,\la)=\left\{\begin{array}{ll}\di\frac{1}{2i}U^0(x,\la)B_1F_1(t\la)
F_2(t\la)U^{0,T}(t,\la)Q(t)&\mbox{for }\;0<t\le x,\\[3mm]
-\di\frac{1}{2i}U^0(x,\la)B_2F^2_1(t\la)U^{0,T}(t,\la)Q(t)&
\mbox{for }\;x<t<a_\la,\\[3mm]
\di\frac{1}{4i}U^0(x,\la)B_2F^2_1(t\la)U^{0,T}(t,\la)BL(t,\la)&
\mbox{for }\;a_\la\le t.
\end{array}\right.
$$
In particular, this yields
$$
\begin{array}{ll}|D_3(x,t,\la)|\le\di\frac{C}{|\la|^{2Re\mu}}t^{-2Re\mu}|Q(t)|&
\mbox{for }\; t<a_\la,\\[3mm]
|D_3(x,t,\la)|\le C|L(t,\la)|&\mbox{for }\; t\ge a_\la.
\end{array}
$$

b) Let $x\ge a_\la.$ Then
$$
D_3(x,t,\la)=\left\{\begin{array}{ll}\di\frac{1}{2i}U^0(x,\la)B_1F_1(t\la)
F_2(t\la)U^{0,T}(t,\la)Q(t) &\mbox{for }\;0<t<a_\la,\\[3mm]
-\di\frac{1}{4i}U^0(x,\la)B_1F_1(t\la)F_2(t\la)U^{0,T}(t,\la)BL(t,\la)&
\mbox{for }\;a_\la\le t<x,\\[3mm]
\di\frac{1}{4i}U^0(x,\la)B_2F^2_1(t\la)\frac{F_2(x\la)}{F_1(x\la)}
U^{0,T}(t,\la)BL(t,\la)& \mbox{for }\;x\le t,
\end{array}\right.
$$
and consequently,
$$
D_3(t,\la)=\left\{
\begin{array}{ll} C|\la t|^{-2Re\mu}|Q(t)| &\mbox{for }\; t<a_\la,\\[3mm]
C|L(t,\la)|&\mbox{for }\; t\ge a_\la.
\end{array}\right.
$$

Using (28) and (29), we calculate
$$
|(U_1)_k(x,\la)|\le C^{k+1}\Big(\frac{1}{|\la|}+\frac{1}{|\la|^{2Re\mu}}
\int_0^{a_\la}t^{-2Re\mu}|Q(t)|\,dt+\int_{a_\la}^\iy|L(t,\la)|\,dt\Big)^k.
$$
Taking lemma 3 into account, we deduce
$$
|(U_1)_k(x,\la)|\le CC^k\Big(\frac{1}{|\la|}+\frac{1}{|\la|}\int_0^\iy
\Big(|Q'(t)|+|Q(t)|\Big)\,dt+\frac{1}{|\la|^\nu} \int_0^\iy
t^{-\nu}|Q(t)|\,dt\Big)^k.
$$
For sufficiently large $|\la|\ge\la_0$, one has $|(U_1)_k(x,\la)|\le C/2^k$.
Therefore, the series $U_1(x,\la)=\di\sum_{k=0}^\iy(U_1)_k(x,\la)$
converges uniformly, hence $U_1(x,\la)$ is continuous with respect to$x,$
and analytic with respect to $\la,$ and $|U_1(x,\la)|\le M_0$.
It follows from (24) and (26) that
$$
|U_1(x,\la)-U^0_1(x,\la)|\le M_0\Big(D_1(\la)+D_2(\la)+\int_0^\iy D_3(t,\la)\,dt\Big).
$$
The theorem is proved.

\medskip
{\bf 4. Asymptotics of the Stockes multipliers. }
Since $E(x,\la)$ and $S(x,\la)$ are fundamental matrices of system (1),
it follows that $E(x,\la)=S(x,\la)\ga(\la)$ and $S(x,\la)=E(x,\la)\be(\la)$;
the matrices $\ga(\la)$ and $\be(\la)$ are called the Stockes multipliers.

\smallskip
{\bf Theorem 6. }{\it The following relations hold:\\
1) $\ga_{j2}(\la)=\la^{\mu_j}\ga^0_{j2},\;j=1,2,$\\
2) $\ga_{j1}(\la)=\la^{\mu_j}\ga^0_{j1}(1+O(|\la|^{-\nu}))$
for $|\la|\to\iy$, $j=1,2,$\\
where $\ga^0_{ij}$ are the Stockes multipliers from $e(x)=C(x)\ga^0$.

Proof. } We rewrite the relations $e(x,\la)=C(x,\la)\ga^0(\la)$ and
$E(x,\la)=S(x,\la)\ga(\la)$ in the vector form:
$$
e_j(x,\la)=\ga^0_{1j}\la^{-\mu}C_1(x,\la)+\ga^0_{2j}\la^\mu C_2(x,\la),
$$
$$
E_j(x,\la)=\ga_{1j}(\la)S_1(x,\la)+\ga_{2j}(\la)S_2(x,\la).
$$
We consider the case $x<a_\la.$ Then $F_j(x\la)=(x\la)^{-\mu}$,
and the last relations imply
$$
\left.\begin{array}{c} U^0_j(x,\la)=\ga^0_{1j}\widehat
C_1(x,\la)+\ga^0_{2j}\cdot(x\la)^{2\mu}
\widehat C_2(x,\la),\\[3mm]
U_j(x,\la)=\ga_{1j}(\la)\la^\mu\widehat S_1(x,\la)+\ga_{2j}(\la)\la^\mu
x^{2\mu}\widehat S_2(x,\la).
\end{array}\right\}                                                    \eqno(30)
$$
Subtracting the first equality from the second one and adding
$\ga^0_{1j}\widehat S_1(x,\la)
-\ga^0_{1j}\widehat S_1(x,\la),$ \\ $\ga^0_{2j}\cdot(x\la)^{2\mu}\widehat
S_2(x,\la)- \ga^0_{2j}\cdot(x\la)^{2\mu}\widehat S_2(x,\la),$ we obtain
$$
U_j(x,\la)-U^0_j(x,\la)=\Big(\ga_{1j}(\la)\la^\mu-\ga^0_{1j}\Big)\widehat
S_1(x,\la) +\ga^0_{1j}\Big(\widehat S_1(x,\la)-\widehat C_1(x,\la)\Big)
$$
$$
+\Big(\ga_{2j}(\la)\la^\mu-\ga^0_{2j}\la^{2\mu}\Big)x^{2\mu}
\widehat S_2(x,\la)+\ga^0_{2j}(x\la)^{2\mu}
\Big(\widehat S_2(x,\la)-\widehat C_2(x,\la)\Big).                    \eqno(31)
$$
For $x\to+0,$ we calculate
$$
U_j(0,\la)-U^0_j(0,\la)=
\Big(\ga_{1j}(\la)\la^\mu- \ga^0_{1j}\Big)\widehat S_1(0,\la).        \eqno(32)
$$
Using (31) we calculate
$$
\Big(\ga_{2j}(\la)\la^\mu-\ga^0_{2j}\la^{2\mu}\Big)\widehat S_2(x,\la)=
\frac{1}{x^{2\mu}}\Big(\Big(U_j(x,\la)-U^0_j(x,\la)\Big)-
\Big(\ga_{1j}(\la)\la^\mu-\ga^0_{1j}\Big)\widehat C_1(x,\la)\Big)
$$
$$
-\frac{1}{x^{2\mu}}\ga_{1j}(\la)\la^\mu\Big(\widehat S_1(x,\la)-
\widehat C_1(x,\la)\Big) -\ga_{2j}(\la)\la^{2\mu}
\Big(\widehat S_2(x,\la)-\widehat C_2(x,\la)\Big).
$$
Taking the estimate $|\widehat S_1(x,\la)-\widehat C_1(x,\la)|\le
Cx^{2Re\mu}\di\int_0^xt^{-2Re\mu}|P(t)|\,dt$ into account, we obtain
$$
(\ga_{2j}(\la)\la^\mu-\ga^0_{2j}\la^{2\mu})\widehat S_2(0,\la)=
\lim_{x\to+0}\frac{1}{x^{2\mu}}\Big((U_j(x,\la)-U^0_j(x,\la))-
(\ga_{1j}(\la)\la^\mu-\ga^0_{1j})\widehat C_1(x,\la)\Big)            \eqno(33)
$$
Since $\widehat S(0,\la)=\widehat C(0,\la)$, $U_{1j}(0,\la)=
U^0_{1j}(0,\la)=0,\;j=1,2,$ it follows from (32)-(33) that
$$
\ga_{1j}(\la)\la^\mu-\ga^0_{1j}=-\frac{1}{c_{10}}
\Big(U_{2j}(0,\la)-U^0_{2j}(0,\la)\Big),                             \eqno(34)
$$
$$
\ga_{2j}(\la)\la^\mu-\ga^0_{2j}\la^{2\mu}=
\lim_{x\to+0}\frac{1}{x^{2\mu}c_{20}}\Big((U_{1j}(x,\la)-U^0_{1j}(x,\la))
-(\ga_{1j}(\la)\la^\mu-\ga^0_{1j})\widehat C_{11}(x,\la)\Big).        \eqno(35)
$$

Let $j=2.$ It follows from (30) that $U_{22}(0,\la)=U^0_{22}(0,\la),$
and, according to (34), $\ga_{12}(\la)\la^\mu-\ga^0_{12}=0.$ Substitute
into (35):
$$
\ga_{22}(\la)\la^\mu-\ga^0_{22}\la^{2\mu}=\lim_{x\to+0}
\Big(\frac{1}{x^{2\mu}c_{20}}\Big(U_{1j}(x,\la)-U^0_{1j}(x,\la)\Big)\Big).
$$
Using
$$
U_2(x,\la)=U^0_2(x,\la)+\di\int_0^x e(x,\la)e^{-1}(t,\la)
\Big(\frac{t}{x}\Big)^{-\mu}BQ(t)U_2(x,\la)\,dt,\; x<a_\la,
$$
and $e(x,\la)=C(x,\la)\ga^0(\la),$ we obtain the estimate
$$
|U_2(x,\la)-U^0_2(x,\la)|\le Cx^{2Re\mu}\di\int_0^xt^{-2Re\mu}|Q(t)|\,dt,
$$
and consequently, $\ga_{22}(\la)\la^\mu-\ga^0_{22}\la^{2\mu}=0.$

Let $j=1.$ For $x<a_\la,$ the equation for $U_1(x,\la)$ has the form
$$
U_1(x,\la)=U^0_1(x,\la)+D_2(x,\la)U_1(a_\la,\la)+
\int_0^{+\iy}D_3(x,t,\la)U_1(t,\la)\,dt.
$$
Taking (34) into account, we calculate $|\ga_{1j}(\la)\la^\mu-
\ga^0_{1j}|\le C|\la|^{-\nu}.$ By virtue of (19), we have
$$
U_{j1}(x,\la)=U^0_{j1}(x,\la)+\Big(e_{j1}(x,\la),\;e_{j2}(x,\la)\Big)
I_1 \int_0^xe^{-1}(t,\la)BQ(t)\Big(\frac{t}{x}\Big)^{-\mu}U_1(t,\la)\,dt
$$
$$
+\frac{1}{2i}\Big(U^0_{j1}(x,\la),\;U^0_{j2}(x,\la)\Big)B_2\Big(
-\int_x^{a_\la}F^2_1(t\la)U^{0,T}(t,\la)Q(t)U_1(t,\la)\,dt
$$
$$
+\int_{a_\la}^\iy F^2_1(t\la)U^{0,T}(t,\la)BL(t,\la)U_1(t,\la)\,dt
+\frac{1}{2}F^2_1(a_\la\la)U^{0,T}(a_\la,\la)BQ^{-1}(a_\la,\la)Q(a_\la)
U_1(a_\la,\la)\Big).
$$
Substituting (34) into (35), we infer
$$
\ga_{21}(\la)\la^\mu-\ga^0_{21}\la^{2\mu}=
\lim_{x\to+0}\frac{1}{x^{2\mu}c_{20}}
\Big(\Big(U_{11}(x,\la)-U^0_{11}(x,\la)\Big)+\frac{1}{c_{10}}
\Big(U_{21}(0,\la)-U^0_{21}(0,\la)\Big)\widehat C_{11}(x,\la)\Big).
$$
Denote
$$
V(\la)=\frac{1}{2i}B_2\Big(
-\int_0^{a_\la}(t\la)^{-2\mu}U^{0,T}(t,\la)Q(t)U_1(t,\la)\,dt+
$$
$$
+\frac{1}{2}\int_{a_\la}^\iy e^{2i\la t}U^{0,T}(t,\la)BL(t,\la)
U_1(t,\la)\,dt+\frac{1}{2}e^{2i\la a_\la}U^{0,T}(a_\la,\la)
BQ^{-1}(a_\la,\la)Q(a_\la)U_1(a_\la,\la)\Big).
$$
Then
$$
\Big(U_{11}(x,\la)-U^0_{11}(x,\la)\Big)+\frac{1}{c_{10}}
\Big(U_{21}(0,\la)-U^0_{21}(0,\la)\Big)\widehat C_{11}(x,\la)
$$
$$
=\Big(e_{11}(x,\la),\;e_{12}(x,\la)\Big)\int_0^x e^{-1}(t,\la)
\Big(\frac{t}{x}\Big)^{-\mu}BQ(t)U_1(t,\la)\,dt
$$
$$
+\Big(U^0_{11}(x,\la)+\frac{U^0_{21}(0,\la)}{c_{10}}\widehat C_{11}(x,\la),
\;\; U^0_{12}(x,\la)+\frac{U^0_{22}(0,\la)}{c_{10}}
\widehat C_{11}(x,\la)\Big)V(\la).
$$
Since $e(x,t)\di\int_0^xe^{-1}(t,\la)\Big(\frac{t}{x}\Big)^{-\mu}BQ(t)
U_1(t,\la)\,dt=\int_0^xG^\br{1}(x,t,\la)BQ(t)U_1(t,\la)\,dt,$ it follows that
$$
\Big|\Big(e_{11}(x,\la),\;e_{12}(x,\la)\Big)\int_0^xe^{-1}(t,\la)
\Big(\frac{t}{x}\Big)^{-\mu}BQ(t)U_1(t,\la)\,dt\Big|\le Cx^{2Re\mu}
\di\int_0^xt^{-2Re\mu}|Q(t)|\,dt.
$$
Furthermore, it follows from (30) that $U_{2j}(0,\la)=-c_{10}\ga^0_{1j}$.
Then
$$
U^0_{1j}(x,\la)+(c_{10})^{-1}U^0_{2j}(0,\la)\widehat C_{11}(x,\la)=
U^0_{1j}(x,\la)-\ga^0_{1j} \widehat C_{11}(x,\la),
$$
and consequently,
$$
U^0_{1j}(x,\la)+(c_{10})^{-1}U^0_{2j}(0,\la)\widehat C_{11}(x,\la)
=\ga^0_{2j}\cdot (x\la)^{2\mu}\widehat C_{12}(x,\la).
$$
It is easy to see that $|V(\la)|\le C|\la|^{-\nu}.$ Thus, we have
$$
\Big|\Big(U_{11}(x,\la)-U^0_{11}(x,\la)\Big)+\frac{1}{c_{10}}
\Big(U_{21}(0,\la)-U^0_{21}(0,\la)\Big)\widehat C_{11}(x,\la)\Big|
$$
$$
\le Cx^{2Re\mu}\Big(\int_0^xt^{-2Re\mu}|P(t)|\,dt+|\la^{2\mu}|\cdot
\frac{1}{|\la|^\nu}\Big),
$$
therefore, $|\ga_{21}(\la)\la^\mu-\ga^0_{21}\la^{2\mu}|\le C|\la^{2\mu}|
\cdot |\la|^{-\nu}.$ The theorem is proved.

\smallskip
{\bf Corollary. }  $|\be_{kj}(\la)-\be^0_{kj}\cdot\la^{-\mu_j}|
\le C|x\la|^{-\nu},\; k,j=1,2$.

\smallskip
{\bf Remark. } Using the above-obtained results, it is easy to deduce
asymptotics of the fundamental matrix $S(x,\la)$ (see [16] for more details):
$$
S_j(x,\la)=\be^0_j\la^{-\mu_j}e^{2i\pi\mu_jm}\left( e^{-i\la
x}\left[\begin{array}{r}-i\\1\end{array}\right]_0-
(-1)^je^{i\pi\mu_jl}
e^{i\la x}\left[\begin{array}{r}i\\1\end{array}\right]_0\right),\;
j=1,2,\;|x\la|\ge1,
$$
$$
\frac{d}{d\la} S_j(x,\la)=\be^0_jx\la^{-\mu_j}e^{2i\pi\mu_jm}\left(e^{-i\la
x}\left[\begin{array}{r}-1\\-i\end{array}\right]_0- (-1)^je^{i\pi\mu_jl}
e^{i\la x}\left[\begin{array}{r}-1\\i\end{array}\right]_0\right),\; |x\la|\ge1,
$$
where\\
$l=\left\{ \begin{array}{rl}1,&\arg(x\la)\in(-\pi,-\pi/2]\cup(\pi/2,\pi],\\
-1,&\arg(x\la)\in(-\pi/2,\pi/2],\end{array}\right.\;\;
m=\left\{ \begin{array}{rl}1,&x<0,\,\arg\la\in(\pi/2,\pi],\\
-1,&x>0,\,\arg\la\in(-\pi,-\pi/2], \\
0, & \mbox{ otherwise,}
\end{array}\right.\;$\\
$\be^0_1\be^0_2=(4i\cos\pi\mu)^{-1}$.

\bigskip
{\bf Acknowledgment.} This work was supported by Grant 1.1436.2014K of the Russian
Ministry of Education and Science and by Grant 13-01-00134 of Russian Foundation for
Basic Research.

\begin{center}
{\bf REFERENCES}
\end{center}
\begin{enumerate}
\item[{[1]}] M.A.Naimark, {\it Linear Differential Operators.} 2nd ed., Nauka,
     Moscow, 1969; English transl. of 1st ed., Parts I,II, Ungar, New York,
     1967, 1968.
\item[{[2]}] Yurko V A, Method of Spectral Mappings in the Inverse
     Problem Theory, Inverse and Ill-posed Problems Series, VSP, Utrecht,
     2002.
\item[{[3]}] Meschanov V.P. and Feldstein A.L., Automatic Design of
     Directional Couplers, Moscow: Sviaz, 1980 (in Russian).
\item[{[4]}] Litvinenko O.N. and Soshnikov V.I., The Theory of
     Heterogenious Lines and their Applications in Radio Engineering,
     Moscow: Radio, 1964 (in Russian).
\item[{[5]}] Freiling G. and Yurko V.A., Reconstructing parameters
     of a medium from incomplete spectral information. Results in
     Mathematics 35 (1999), 228-249.
\item[{[6]}] Anderssen R.S., The effect of discontinuities in density
     and shear velocity on the asymptotic overtone structure of
     tortional eigenfrequencies of the Earth. Geophys. J.R. astr. Soc.
     50 (1997), 303-309.
\item[{[7]}] Lapwood F.R. and Usami T., Free Oscilations of the Earth,
     Cambridge University Press, Cambridge, 1981.
\item[{[8]}] Hald O.H., Discontinuous inverse eigenvalue problems.
     Comm. Pure Appl. Math. 37 (1984), 539-577.
\item[{[9]}] Constantin A., On the inverse spectral problem for the
     Camassa-Holm equation. J. Funct. Anal. 155 (1998), no. 2, 352-363.
\item[{[10]}] Gasymov M.G. Determination of Sturm-Liouville equation with
     a singular point from two spectra. Doklady Akad. Nauk SSSR 161 (1965),
     274-276; transl. in Sov. Math. Dokl. 6(1965), 396-399.
\item[{[11]}] Zhornitskaya L.A. and Serov V.S., Inverse eigenvalue
     problems  for  a singular  Sturm-Liouville  operator  on  (0,1).
     Inverse Problems 10 (1994), no.4, 975-987.
\item[{[12]}]  Yurko V.A., Inverse problem for differential equations
     with a singularity. Differen. Uravneniya 28 (1992), 1355-1362;
     English transl. in Differential Equations  28 (1992), 1100-1107.
\item[{[13]}] Yurko V.A., On higher-order differential operators with a
     singular point. Inverse Problems 9 (1993), 495-502.
\item[{[14]}] Yurko V.A., On higher-order differential operators with a
     regular  singularity.  Mat. Sb. 186 (1995),  no.6,  133-160;
     English  transl. in  Sbornik; Mathematics  186 (1995), no.6, 901-928.
\item[{[15]}] Yurko V.A., Integral transforms connected with
     differential operators having singularities inside the interval.
     Integral Transforms and Special Functions 5 (1997), no.3-4, 309-322.
\item[{[16]}] Gorbunov O. Inverse problem for Dirac operators with 
     non-integrable singularities inside the interval. PhD Thesis. 
     Saratov University, Saratov, Russia, 2003.
\end{enumerate}

\begin{tabular}{ll}
Name:             &   Yurko, Vjacheslav  \\
Place of work:    &   Department of Mathematics, Saratov State University \\
{}                &   Astrakhanskaya 83, Saratov 410012, Russia \\
E-mail:           &   yurkova@info.sgu.ru  \\
\end{tabular}

\begin{tabular}{ll}
Name:             &   Gorbunov, Oleg  \\
Place of work:    &   Department of Mathematics, Saratov State University \\
{}                &   Astrakhanskaya 83, Saratov 410012, Russia \\
E-mail:           &   gorbunovob@info.sgu.ru  \\
\end{tabular}

\begin{tabular}{ll}
Name:             &   Shieh, Chung-Tsun  \\
Place of work:    &   Tamkang University \\
{}                &   Taiwan \\
E-mail:           &   ctshieh@mail.tku.edu.tw   \\
\end{tabular}

\end{document}